%&plain

%%%%%%%   Text for my Bourbaki talk, given 12 June 1999
%%%%%%%    --- Robert L. Bryant   Today:  25 August 1999

\magnification 1200
\topskip=20pt
\vsize=19,5cm
\hsize=13cm

\openup 1\jot

\def\expno{861}

\headline{\ifnum\pageno=1 \hbox{}\hfill\else
\ifnum\pageno<10 \hbox to\hsize{\tenrm\hfil
\expno-0\folio \hfil}\else
\hbox to\hsize{\tenrm\hfil \expno-\folio \hfil}\fi\fi}
\nopagenumbers

% I need some of the standard AMS Fonts and symbols:

\input amssym.def
\input amssym.tex

%definition of small capitals for the bibliography
\newfam\smcfam%
\font\tensmc=cmcsc10%
\def\smc{\fam\smcfam\tensmc}%\smc is family 7
\textfont\smcfam=\tensmc%

%definition of the font for footnotes
\font\ninerm=cmr9%
\font\ninebf=cmbx9%
%

% font for Schur functor symbols: CM sans serif
\font\sans=cmss10
\def\Lam{\hbox{\sans\char3}}%sans Lambda
\def\Sym{\hbox{\sans S}}

\def\ts{\textstyle }

\def\bbR{{\Bbb R}}
\def\bbC{{\Bbb C}}
\def\bbF{{\Bbb F}}
\def\bbH{{\Bbb H}}

\def\bbP{{\Bbb P}}
\def\bbZ{{\Bbb Z}}

\def\g{\gamma}

\def\defeq{\mathbin{\,\raise0.29pt\hbox{:}\hskip-1.0pt{=}\,}}

\def\w{{\mathchoice{\,{\scriptstyle\wedge}\,}
{{\scriptstyle\wedge}}
{{\scriptscriptstyle\wedge}}{{\scriptscriptstyle\wedge}}}}

\def\eug{{\frak g}}\def\euh{{\frak h}}
\def\eus{{\frak s}}\def\eum{{\frak m}}
\def\eugl{{\frak{gl}}}
\def\eusl{{\frak{sl}}}\def\euso{{\frak{so}}}

\def\SO{{\rm SO}}\def\SL{{\rm SL}}
\def\Sp{{\rm Sp}}\def\SU{{\rm SU}}
\def\CO{{\rm CO}}\def\GL{{\rm GL}}
\def\Spin{{\rm Spin}}\def\U{{\rm U}}
\def\E{{\rm E}}\def\O{{\rm O}}
\def\Diff{{\rm Diff}}\def\G{{\rm G}}

\def\Ad{{\rm Ad}}\def\Hom{{\rm Hom}}
\def\coker{{\rm coker}}

\def\pbold{{\bf p}}

%%%%%  This is the end of the macros.  Now the paper begins.

\noindent
\line{S\'eminaire BOURBAKI \hfill Juin 1999}\par
\noindent 51\`eme ann\'ee, 1998-99, n$^{\rm o}$ \expno
\vskip 1.9truecm

\centerline{\bf RECENT ADVANCES IN THE THEORY OF HOLONOMY} 
\medskip
\centerline {by {\ninebf Robert BRYANT}}
\vskip 1.9truecm

\noindent
{\bf 1. I\ninebf NTRODUCTION}

\vskip 0.5truecm
\noindent
{\bf 1.0.} {\smc Outline.} This report is organized as follows:

\itemitem{1.} Introduction
\itemitem{2.} Riemannian Holonomy
\itemitem{3.} Torsion-free non-metric connections: the irreducible case

In a short lecture of this nature, it is impossible to describe the
history of the subject in any depth, but the reader can find more
information on the Riemannian (and pseudo-Riemannian) case by 
consulting [Bes], [Sa], and the forthcoming, much-anticipated~[Jo3], 
especially for the exceptional cases.  For the non-metric case,
aside from the survey~[Br3], the expository papers~[MS2] and [Schw1]
provide a valuable account of both the representation theoretic
and twistor theoretic approaches to the study of holonomy.

\vskip 0.5truecm
\noindent
{\bf 1.1.} {\smc Historical Remarks.}  According to the Oxford English
Dictionary, it was Heinrich Hertz in~1899 who introduced the words 
{\it holonomic\/} and {\it nonholonomic\/} to describe a property of 
velocity constraints in mechanical systems.  

Velocity constraints are {\it holonomic\/} if they force a curve in 
state space to stay in a proper subspace.  As an example, the 
condition~$\pbold\cdot d\pbold = 0$ for a vector particle~$\pbold\in\bbR^n$
forces~$\pbold$ to have constant length, while the 
constraint~$\pbold\wedge d\pbold = 0$ forces~$\pbold$ to move on a line.

{\it Nonholonomic\/} constraints, on the other hand, imply no such 
`finite' constraints.  A classical example is that of a ball rolling
on a table without slipping or twisting.  The state space
is~$B=\SO(3)\times\bbR^2$, where the~$\SO(3)$
records the orientation of the ball and the~$\bbR^2$
records its contact point on the plane.  The rolling constraint
is expressed as the set of differential equations
$$
\alpha \defeq a^{-1}\,da+a^{-1}\pmatrix{0&0&-dx\cr0&0&-dy\cr dx&dy&0\cr}a = 0
\leqno{(0.1)}
$$
for a curve~$\bigl(a(t);x(t),y(t)\bigr)$ in~$B$.  The curves in~$B$ satisfying
this constraint are those tangent to the $2$-plane field~$D=\ker\alpha$ 
that is transverse to the fibers of the 
projection~$\SO(3)\times\bbR^2\to\bbR^2$. It is not difficult to show that 
any two points of~$B$ can be joined by a curve tangent to~$D$.  

Such constraints and their geometry have long been of considerable
interest in the calculus of variations and control theory.  
For recent results, see the foundational work on Carnot-Carath\'eodory
geometries by Gromov~[Gr]. 

\vskip 0.5truecm
\goodbreak
\noindent
{\bf 1.2.} {\smc The holonomy group.}
\'Elie Cartan~[Ca2] introduced the {\it holonomy group\/} in the context 
of differential geometry. It measures the failure of the parallel translation 
associated to a connection to be holonomic.  

The data are a connected manifold~$M$,
a Lie group~$H$ with Lie algebra~$\euh$, a principal right 
$H$-bundle~$\pi:B\to M$, and a connection~$\theta$ on~$B$, i.e., $\theta$ is 
an $\euh$-valued \hbox{$1$-form} on~$B$ that pulls back to each 
$\pi$-fiber to be the canonical left invariant 1-form on~$H$ and that
satisfies the equivariance relation~$R_h^*(\theta) = \Ad(h^{-1})(\theta)$. 
(The nonholonomic example described above with $(M,H,B,\theta)
=\bigl(\bbR^2,\,\SO(3),\,\SO(3){\times}\bbR^2,\,\alpha\bigr)$ is
an example.)

A piecewise $C^1$ curve~$\g:[0,1]\to B$ is said to be 
{\it $\theta$-horizontal} or {\it $\theta$-parallel} if~$\g^*(\theta)=0$.
The {\it $\theta$-holonomy\/}~$B^\theta_u\subset B$ of~$u\in B$ 
is defined to be the set of values of~$\g(1)$ as~$\g:[0,1]\to B$
ranges over the $\theta$-horizontal curves with~$\g(0)=u$.  
The $\Ad(H)$-equivariance of~$\theta$ (coupled with the connectedness
of~$M$) implies that there is a 
subgroup~$H^\theta_u\subset H$ so that~$B^\theta_u$ is an 
$H^\theta_u$-subbundle of~$B$ and that $H^\theta_{u\cdot h} 
= h^{-1}\,H^\theta_u\,h$ for~$h\in H$. Consequently, the 
conjugacy class of~$H^\theta_u\subset H$ depends only on~$\theta$.  The 
group~$H^\theta_u$ (or, more informally, its conjugacy class in~$H$) 
is called the {\it holonomy\/} of~$\theta$.

It is a theorem of Borel and Lichnerowicz~[BL] that~$H^\theta_u$ is 
a Lie subgroup of~$H$.  By a theorem of Ambrose 
and Singer~[AS], the Lie algebra~$\euh^\theta_u$ 
of~$H^\theta_u$ is spanned by the set
$$
\{\ \Theta(x,y)\ \vrule\ \ x,y\in T_vB\,,\ v\in B^\theta_u\ \}
$$
where~$\Theta = d\theta + {1\over2}[\theta,\theta]$ is the {\it curvature\/}
of~$\theta$.  The identity component $(H^\theta_u)^0\subset H^\theta_u$ 
is known as the {\it restricted holonomy\/} of~$\theta$.  There is a 
well-defined surjective homomorphism~$\rho^\theta:\pi_1\bigl(M,\pi(u)\bigr)
\to H^\theta_u/(H^\theta_u)^0$ that 
satisfies~$\rho^\theta\bigl([\pi{\circ}\g]\bigr) = \g(1)(H^\theta_u)^0$
for every $\theta$-horizontal curve~$\g$ with~$\g(0)=u$ 
and~$\g(1)\in u{\cdot}H$.

Using these results, it can be shown~[KNo] that, when~$\dim M>1$, 
a Lie subgroup~$G\subset H$ can be the holonomy group of a connection 
on~$B$ if and only if~$B$ admits a structural reduction to a $G$-bundle.
Thus, the set of possible holonomies of connections on~$B$ is determined
topologically.

\vskip 0.5truecm
\noindent
{\bf 1.3.} {\smc $H$-structures and torsion.}  A common source of
connections in geometry is that of~$H$-structures on manifolds.
Suppose $\dim M = n$ and let~$\eum$ be a reference vector space of 
dimension~$n$.  An ($\eum$-valued) {\it coframe\/} at~$x\in M$ is a linear
isomorphism~$u:T_xM\to \eum$.  The set~$F(M,\eum)$ of $\eum$-valued coframes 
at the points of~$M$ is naturally a principal right $\GL(\eum)$-bundle 
over~$M$, with basepoint projection~$\pi:F(M,\eum)\to M$. There is a 
tautological $\eum$-valued 1-form~$\omega$ on~$F(M,\eum)$ defined by the 
formula~$\omega(v) = u\bigl(\pi'(v)\bigr)$ for all~$v\in T_uF(M,\eum)$.

Let~$H\subset\GL(\eum)$ be a subgroup.  An {\it $H$-structure\/} on~$M$ 
is an $H$-subbundle~$B$ such that~$B\subset F(M,\eum)$.  
When~$H$ is a closed subgroup 
of~$\GL(\eum)$ (the only case I will consider today), the set 
of~$H$-structures on~$M$ is the set of sections of the 
bundle~$F(M,\eum)/H\to M$.  The problem of determining whether 
there exists an $H$-structure on~$M$ is a purely topological one. 

Most of the familiar geometric
structures on~$M$ can be described as $H$-structures.  
For example, when $H=\O(Q)\subset\GL(\eum)$ is the group of linear 
transformations preserving a quadratic form~$Q$ of type~$(p,q)$ on~$\eum$, 
a choice of $H$-structure on~$M$ is equivalent to a choice of 
pseudo-Riemannian metric of type~$(p,q)$ on~$M$.  
When~$H=\Sp(S)\subset\GL(\eum)$ is the group of linear 
transformations preserving a nondegenerate skewsymmetric bilinear form~$S$ 
on~$\eum$, a choice of $H$-structure on~$M$ is equivalent to a choice 
of a nondegenerate $2$-form~$\sigma$ on~$M$, i.e., an almost symplectic 
structure. 

If~$\pi:B\to M$ is an $H$-structure on~$M$, the tautological form~$\omega$
pulled back to~$B$ will also be denoted~$\omega$ when there is no chance
of confusion.  If~$\theta$ is a connection on~$B$, then the {\it first
structure equation of Cartan\/} says that there exists an $H$-equivariant 
function~$T:B\to\Hom\bigl(\Lam^2(\eum),\eum\bigr)\simeq
\eum\otimes\Lam^2(\eum^*)$ so that
$$
d\omega + \theta\w\omega = {\ts{1\over2}}\,T(\omega\w\omega).
\leqno{(0.2)}
$$
The function~$T$ is the {\it torsion function\/} of~$\theta$,
and $\theta$ is {\it torsion-free\/} if~$T=0$.  Any 
connection~$\theta'$ on~$B$ is~$\theta + p(\omega)$
where~$p:B\to\Hom(\eum,\euh)$ is an $H$-equivariant function. 
Its torsion function is~$T' = T + \delta(p)$, where~$\delta:\euh\otimes\eum^*
\to\eum{\otimes}\Lam^2(\eum^*)$ is defined via the inclusion~$\euh\subset
\eum\otimes\eum^*$ and the skewsymmetrizing 
map~$\eum{\otimes}\eum^*{\otimes}\eum^*\to\eum{\otimes}\Lam^2(\eum^*)$.

Evidently, the reduced map~$\bar T:B\to \coker(\delta)$ is independent of
the choice of~$\theta$ and is the obstruction to choosing a torsion-free
connection on~$B$.  For example, when~$H=\Sp(S)$ as above, it is not
difficult to see that~$\coker(\delta)\simeq\Lam^3(\eum^*)$ and that the
reduced torsion represents the exterior derivative~$d\sigma$ of the 
associated almost symplectic form~$\sigma\in\Omega^2(M)$. 
Thus, an~$\Sp(S)$-structure admits a torsion-free
connection if and only if the almost symplectic structure on~$M$ is actually
symplectic.  Furthermore, $\ker(\delta)\simeq\Sym^3(\eum^*)$, so that 
the torsion function~$T$ does not determine a unique connection~$\theta$.
By contrast, when~$H=\O(Q)$, the map~$\delta$ is an 
isomorphism, which is merely the fundamental lemma of Riemannian geometry:  
Every pseudo-Riemannian metric on~$M$ possesses a unique 
compatible, torsion-free connection.

The reduced torsion is the first order local invariant for
$H$-structures on~$M$ and $H$-structures satisfying~$\bar T=0$ are usually
referred to as {\it torsion-free\/} or {\it $1$-flat\/}.  Many (but not all)
of the $H$-structures that have received the most attention in differential
geometry are $1$-flat.  All pseudo-Riemannian metric structures
are $1$-flat, an almost symplectic structure is $1$-flat if and only if 
it is symplectic (Darboux' Theorem), an almost complex structure is $1$-flat 
if and only if it is complex (Newlander-Nirenberg Theorem), an Hermitian
structure is $1$-flat if and only if it is K\"ahler.  Contact structures
and Carnot-Carath\'eodory structures, on the other hand, are definitely
not $1$-flat, as the nondegeneracy of the reduced torsion is an essential 
part of the geometry.
 
\vskip 0.5truecm
\noindent
{\bf 1.4.} {\smc Criteria for holonomy in the torsion-free case.}
The condition~$\bar T=0$ is a $\Diff(M)$-invariant, 
first order equation for $H$-structures on~$M$.  
Given a torsion-free $H$-structure~$B\subset F(M,\eum)$, 
one can ask about the possibilities for the holonomy group of
a torsion-free connection~$\theta$ on~$B$.  When $\ker(\delta)=0$, a
very common situation, the torsion-free connection, if it exists, 
will be unique, so that it makes sense to speak of the holonomy of~$B$ itself.

The vanishing of the torsion implies nontrivial restrictions on the
holonomy.  Assuming~$T=0$
and differentiating~(0.2) yields $\Theta\w\omega=0$.  Thus, the
{\it curvature function\/}~$R:B\to\euh\otimes\Lam^2(\eum^*)$ of~$\theta$, 
for which the {\it second structure equation of Cartan}
$$
\Theta = d\theta + {\ts{1\over2}}\,[\theta,\theta] 
    = {\ts{1\over2}} R(\omega\w\omega)
\leqno{(0.3)}
$$
holds, takes values in the kernel~$K(\euh)$ of the
map~$\delta:\euh\otimes\Lam^2(\eum^*)\to \eum\otimes\Lam^3(\eum^*)$ 
defined by the same methods%
\footnote{$^1$}{\baselineskip=9pt{\ninerm The
maps I am denoting by~$\delta$ (as well as the ones to follow)
are part of the {\it Spencer complex\/}
associated to the inclusion~$\euh\subset\eum\otimes \eum^*$, see~[Br3].
These maps are~$\euh$-module maps and so the various kernels and
cokernels to be introduced are $\euh$-modules as well.}\par}
 as the previous~$\delta$.  In particular, if there is
a proper subalgebra~$\eug\subset\euh$ so 
that~$K(\euh)\subseteq\eug\otimes\Lam^2(\eum^*)$, then, by the Ambrose-Singer
holonomy theorem, the (restricted) holonomy of~$\theta$ must lie in 
the connected Lie group~$G\subset H$ whose Lie algebra is~$\eug$.  

The intersection of the subspaces~$\eus\subset\euh$ that satisfy
$K(\euh)\subseteq\eus\otimes\Lam^2(\eum^*)$ is an ideal~$\eug\subset\euh$.
Thus, a necessary condition that there exist
a torsion-free connection with holonomy~$H\subset\GL(\eum)$ is that
$K(\euh)\not=K(\eug)$ for any proper ideal~$\eug\subset\euh$.  
This is usually referred to as {\it Berger's first criterion\/}~[Ber, Br3].

This criterion is very restrictive:  If~$\euh$ is semi-simple,
there are, up to equivalence, only a finite number of 
representations~$\euh\hookrightarrow\eugl(\eum)$ without trivial summands
satisfying it.  As a simple example, if $\euh\simeq\eusl(2,\bbR)$, 
and~$V_k\simeq\Sym^k(V_1)$ denotes the irreducible 
$\eusl(2,\bbR)$-representation of dimension~$k{+}1$, then the
$\eusl(2,\bbR)$-representations~$\eum$ without~$V_0$-summands that satisfy
Berger's first criterion are~$V_1$,\ $V_1{\oplus}V_1$,\ $V_2$,\ $V_3$,\ 
and~$V_4$.

{\it Symmetric examples.}
One large class of examples where Berger's first criterion is satisfied
is provided by the following construction:  Suppose that there is a 
surjective skewsymmetric pairing~$\eum\times \eum\to\euh$ so that, together 
with the Lie algebra bracket on~$\euh$ and the $\euh$-module pairing
$\euh\times\eum\to\eum$, it defines a Lie algebra on~$\eug=\euh{\oplus}\eum$.  
Then the pair~$(\eug,\euh)$ is a symmetric
pair of Lie algebras~[KN].  Let~$G$ be the simply connected Lie group 
with Lie algebra~$\eug$ and let~$\tilde H\subset G$ be the (necessarily 
closed) connected subgroup corresponding to the subalgebra~$\euh\subset\eug$
and let~$H=\Ad_\eum(\tilde H)\subset\GL(\eum)$ be its almost faithful
image.  Then~$M = G/\tilde H$ is an affine symmetric space in a canonical way,
and the coset projection~$G\to M$ covers a torsion-free $H$-structure on~$M$ 
with connection whose holonomy is~$H$ (by the assumption~$[\eum,\eum]=\euh$ 
and the Ambrose-Singer holonomy theorem.)  The classification of the 
symmetric Lie algebra pairs~$(\eug,\euh)$ is a (rather involved) algebra 
problem.  It was solved by Berger~[Ber2] in the case that~$\eug$ is 
semi-simple or when~$\euh$ acts irreducibly on~$\eum\simeq\eug/\euh$. 

The case where~$\theta$ is a locally symmetric connection with holonomy~$H$
is characterized on the $H$-structure~$B$ by the condition 
that~$R:B\to K(\euh)$ be constant. In fact, differentiating~(0.3)
yields~$0 = (dR{+}\theta{.}R)(\omega\w\omega)$ (where~$\theta{.}R$ is 
the result of the $\euh$-module pairing~$\euh{\times}K(\euh)\to K(\euh)$).
Equivalently, $dR = -\theta{.}R + R'(\omega)$ 
where~$R':B\to K(\euh)\otimes \eum^*$ must take values in the 
kernel~$K^1(\euh)$ of the natural 
map~$\delta:K(\euh)\otimes \eum^*\to\euh\otimes\Lam^3(\eum^*)$. Now,
$R'$ vanishes identically exactly when~$R$ is parallel, which is 
exactly when the pair~$(B,\theta)$ defines a locally symmetric affine
structure on~$M$.  Thus, if~$H$ can be the holonomy of a torsion-free
affine connection that is not locally symmetric, then~$K^1(\euh)\not=0$.  
This is {\it Berger's second criterion}. 

For example, when~$\euh=\eusl(2,\bbR)$ and~$\eum=V_4$, one has~$K^1(\euh)=0$.
Thus, this $5$-dimensional representation of~$\SL(2,\bbR)$ 
could occur as holonomy of a torsion-free connection on~$M^5$ only
when that connection is  locally symmetric.  
In fact, it occurs as the holonomy of the symmetric 
spaces~$\SL(3,\bbR)/\SO(2,1)$ and~$\SU(2,1)/\SO(2,1)$ and in no other way.
The other four cases with~$\euh=\eusl(2,\bbR)$ satisfying Berger's first 
criterion also satisfy Berger's second criterion.  Moreover, each does occur 
as the holonomy of a torsion-free connection that is not locally symmetric.

\vskip 0.5truecm
\noindent
{\bf 1.5.} {\smc Classification.}
In the case where~$\eum$ is an irreducible $\euh$-module 
(implying that $\euh$ is reductive),
Berger's fundamental works~[Ber1, Ber2] went a long way towards
classifying the subalgebras~$\euh\subset\eugl(\eum)$ that satisfy
his first and second criteria.  The methods involved heavy use of
representation theory and ultimately reduced to a painstaking, elaborate 
case analysis.  This list was refined and completed
by the combined work of several people:  Alekseevskii~[Al1],
myself, and most recently and importantly, the combined work of Chi, 
Merkulov, and Schwachh\"ofer.  I will discuss this further in~\S3.

Berger's work provided a (partial) list of candidates for the
irreducibly acting holonomy groups of torsion-free connections
that are not locally symmetric.  
He divided the list into two parts: The first part consists 
of~$\euh\subset\eugl(\eum)$ that lie in some~$\euso(Q)$ for some 
nondegenerate quadratic form~$Q$ on~$\eum$, so that the associated 
$H$-structure defines a pseudo-Riemannian structure on~$M$. For
these cases, the injectivity of the map~$\delta:\euh\otimes\eum^*
\to\eum\otimes\Lam^2(\eum^*)$ implies that the torsion-free connection
is unique, so that it makes sense to speak of the holonomy of the
underlying $H$-structure itself.  I will refer to this part as the 
{\it metric list\/}.  The second part consists of~$\euh\subset\eugl(\eum)$ 
that do not lie in any~$\euso(Q)$ and hence will be referred to as
the {\it nonmetric list\/}.  For many of the
subalgebras on the nonmetric list, the map~$\delta$ is not injective,
so that the associated~$H$-structure does not determine the torsion-free
connection. 

\topinsert
\centerline{\bf I. Pseudo-Riemannian, Irreducible Holonomies in~$\bbR^n$}
\nobreak
\centerline{
\vbox{\offinterlineskip
\halign{
&\vrule#&\strut\quad\hfil#\hfil\quad\cr  %The preamble
\multispan6{\hrulefill}&\cr
height 2 pt&\omit&&\omit&&\omit&\cr
& $n$ && H && Local Generality* &\cr
height 2 pt&\omit&&\omit&&\omit&\cr
\multispan6{\hrulefill}&\cr
height 1 pt&\omit&&\omit&&\omit&\cr
\multispan6{\hrulefill}&\cr
height 3 pt&\omit&&\omit&&\omit&\cr
& $p{+}q\ge2$ && $\SO(p,q)$ && ${1\over2}n(n{-}1)$ of~$n$ &\cr
height 2 pt&\omit&&\omit&&\omit&\cr
& $2p$ && $\SO(p,\bbC)$ && ${1\over2}p(p{-}1)^\bbC$ of~$p^\bbC$ &\cr
height 3 pt&\omit&&\omit&&\omit&\cr
\multispan6{\hrulefill}&\cr
height 3 pt&\omit&&\omit&&\omit&\cr
& $2(p{+}q)\ge4$ && $\U(p,q)$ && 1 of~$n$  &\cr
height 2 pt&\omit&&\omit&&\omit&\cr
& $2(p{+}q)\ge4$ && $\SU(p,q)$ && 2 of~$n{-}1$  &\cr
height 3 pt&\omit&&\omit&&\omit&\cr
\multispan6{\hrulefill}&\cr
height 3 pt&\omit&&\omit&&\omit&\cr
& $4(p{+}q)\ge8$ && $\Sp(p,q)$ && $2(p{+}q)$ of $(2p{+}2q{+}1)$ &\cr
height 3 pt&\omit&&\omit&&\omit&\cr
\multispan6{\hrulefill}&\cr
height 3 pt&\omit&&\omit&&\omit&\cr
& $4(p{+}q)\ge8$ && $\Sp(p,q)\!\cdot\!\Sp(1)$ && $2(p{+}q)$ of 
$(2p{+}2q{+}1)$ &\cr
height 2 pt&\omit&&\omit&&\omit&\cr
& $4p\ge8$ && $\Sp(p,\bbR)\!\cdot\!\SL(2,\bbR)$ && $2p$ of (2p{+}1) &\cr
height 2 pt&\omit&&\omit&&\omit&\cr
& $8p\ge16$ && $\Sp(p,\bbC)\!\cdot\!\SL(2,\bbC)$ && $2p^\bbC$ of 
$(2p{+}1)^\bbC$ &\cr
height 3 pt&\omit&&\omit&&\omit&\cr
\multispan6{\hrulefill}&\cr
height 3 pt&\omit&&\omit&&\omit&\cr
& 7 && $G_2$ && 6 of 6 &\cr
height 2 pt&\omit&&\omit&&\omit&\cr
& 7 && $G'_2$ && 6 of 6 &\cr
height 2 pt&\omit&&\omit&&\omit&\cr
& 14 && $G_2^\bbC$ && $6^\bbC$ of~$6^\bbC$ &\cr
height 3 pt&\omit&&\omit&&\omit&\cr
\multispan6{\hrulefill}&\cr
height 3 pt&\omit&&\omit&&\omit&\cr
& 8 && $\Spin(7)$ && 12 of 7 &\cr
height 2 pt&\omit&&\omit&&\omit&\cr
& 8 && $\Spin(4,3)$ && 12 of 7 &\cr
height 2 pt&\omit&&\omit&&\omit&\cr
& 16 && $\Spin(7,\bbC)$ && $12^\bbC$ of~$7^\bbC$ &\cr
height 3 pt&\omit&&\omit&&\omit&\cr
\multispan6{\hrulefill}&\cr
\noalign{\smallskip}
\omit \vbox{\hsize=4.9truein \baselineskip=12pt
{\ninerm \noindent *Counted modulo diffeomorphism. 
The notation ``$d$ of $\ell$'' means ``$d$ functions of $\ell$ variables''
and a superscript~$\bbC$ means `holomorphic'.}}
\hidewidth\cr
\noalign{\vskip-10pt}
}
}
}
\endinsert

\topinsert
\centerline{\bf II. The `Classical' Non-Metric, Irreducible Holonomies}
\medskip\nobreak
\centerline{
\vbox{\offinterlineskip
\halign{
&\vrule#&\strut\quad\hfil#\hfil\quad\cr  %The preamble
\multispan6{\hrulefill}&\cr
height 2 pt&\omit&&\omit&&\omit&\cr
& H$^\lozenge$ && $\eum$ && Restrictions$^\triangledown$ &\cr
height 2 pt&\omit&&\omit&&\omit&\cr
\multispan6{\hrulefill}&\cr
height 2 pt&\omit&&\omit&&\omit&\cr
\multispan6{\hrulefill}&\cr
height 3 pt&\omit&&\omit&&\omit&\cr
& $G_\bbR\!\cdot\!\SL(n,\bbR)$ && $\bbR^n$ &&$n\ge2$&\cr
height 2 pt&\omit&&\omit&&\omit&\cr
& $G_\bbC\!\cdot\!\SL(n,\bbC)$ && $\bbC^n$ &&$n\ge1$&\cr
height 2 pt&\omit&&\omit&&\omit&\cr
& $G_\bbR\!\cdot\!\SL(n,\bbH)$ && $\bbH^n$ &&$n\ge1$&\cr
height 3 pt&\omit&&\omit&&\omit&\cr
\multispan6{\hrulefill}&\cr
height 3 pt&\omit&&\omit&&\omit&\cr
& $\Sp(n,\bbR)$ && $\bbR^{2n}$ &&$n\ge2$&\cr
height 2 pt&\omit&&\omit&&\omit&\cr
& $\Sp(n,\bbC)$ && $\bbC^{2n}$ &&$n\ge2$&\cr
height 3 pt&\omit&&\omit&&\omit&\cr
\multispan6{\hrulefill}&\cr
height 3 pt&\omit&&\omit&&\omit&\cr
& $\bbR^+\!\cdot\!\Sp(2,\bbR)$ && $\bbR^{4}$ && &\cr
height 2 pt&\omit&&\omit&&\omit&\cr
& $\bbC^*\!\cdot\!\Sp(2,\bbC)$ && $\bbC^{4}$ && &\cr
height 3 pt&\omit&&\omit&&\omit&\cr
\multispan6{\hrulefill}&\cr
height 3 pt&\omit&&\omit&&\omit&\cr
& $\CO(p,q)$ && $\bbR^{p+q}$ &&$p+q\ge3$&\cr
height 2 pt&\omit&&\omit&&\omit&\cr
& $G_\bbC\!\cdot\!\SO(n,\bbC)$ && $\bbC^n$ &&$n\ge3$&\cr
height 3 pt&\omit&&\omit&&\omit&\cr
\multispan6{\hrulefill}&\cr
height 3 pt&\omit&&\omit&&\omit&\cr
& $G_\bbR\!\cdot\!\SL(p,\bbR)\!\cdot\!\SL(q,\bbR)$ &
  & $\bbR^{pq\phantom{2}}\simeq\bbR^p{\otimes}_\bbR\bbR^q$ &
  &$p\ge q\ge2,\ (p,q)\not=(2,2)$&\cr
height 2 pt&\omit&&\omit&&\omit&\cr
& $G_\bbC\!\cdot\!\SL(p,\bbC)\!\cdot\!\SL(q,\bbC)$ &
  & $\bbC^{pq\phantom{2}}\simeq\bbC^p{\otimes}_\bbC\bbC^q$ &
  &$p\ge q\ge2,\ (p,q)\not=(2,2)$&\cr
height 2 pt&\omit&&\omit&&\omit&\cr
& $G_\bbR\!\cdot\!\SL(p,\bbH)\!\cdot\!\SL(q,\bbH)$ &
  & $\bbR^{4pq}\simeq \bbH^p{\otimes}_\bbH\bbH^q$ &
  &$p\ge q\ge1,\ (p,q)\not=(1,1)$&\cr
height 3 pt&\omit&&\omit&&\omit&\cr
\multispan6{\hrulefill}&\cr
height 3 pt&\omit&&\omit&&\omit&\cr
& $G_\bbR\!\cdot\!\SL(p,\bbC)$ &
  & $\bbR^{p^2}\simeq\bigl(\bbC^p{\otimes}_\bbC\overline{\bbC^p}\bigr)^\bbR$ 
  &&$p\ge 3$&\cr
height 3 pt&\omit&&\omit&&\omit&\cr
\multispan6{\hrulefill}&\cr
height 3 pt&\omit&&\omit&&\omit&\cr
& $G_\bbR\!\cdot\!\SL(p,\bbR)$ &
  & $\bbR^{p(p+1)/2}\simeq \Sym^2_\bbR(\bbR^p)$ 
  &&$p\ge 3$&\cr
height 2 pt&\omit&&\omit&&\omit&\cr
& $G_\bbC\!\cdot\!\SL(p,\bbC)$ &
  & $\bbC^{p(p+1)/2}\simeq\Sym^2_\bbC(\bbC^p)$ &
  &$p\ge 3$&\cr
height 2 pt&\omit&&\omit&&\omit&\cr
& $G_\bbR\!\cdot\!\SL(p,\bbH)$ &
  & $\bbR^{p(2p+1)\hphantom{/}}\simeq\Sym^2_\bbH(\bbH^p)$ &
  &$p\ge 2$&\cr
height 3 pt&\omit&&\omit&&\omit&\cr
\multispan6{\hrulefill}&\cr
height 3 pt&\omit&&\omit&&\omit&\cr
& $G_\bbR\!\cdot\!\SL(p,\bbR)$ &
  & $\bbR^{p(p-1)/2}\simeq\Lam^2_\bbR(\bbR^p)$ &
  &$p\ge 5$&\cr
height 2 pt&\omit&&\omit&&\omit&\cr
& $G_\bbC\!\cdot\!\SL(p,\bbC)$ &
  & $\bbC^{p(p-1)/2}\simeq\Lam^2_\bbC(\bbC^p)$ &
  &$p\ge 5$&\cr
height 2 pt&\omit&&\omit&&\omit&\cr
& $G_\bbR\!\cdot\!\SL(p,\bbH)$ &
  & $\bbR^{p(2p-1)\hphantom{/}}\simeq\Lam^2_\bbH(\bbH^p)$ &
  &$p\ge 3$&\cr
height 3 pt&\omit&&\omit&&\omit&\cr
\multispan6{\hrulefill}&\cr
\noalign{\smallskip}
\omit\ 
{\ninerm $^\lozenge$ $G_\bbF$ is any connected subgroup of~$\bbF^*$.
\qquad\qquad
 $^\triangledown$ To avoid repetition or reducibility.}
\hidewidth\cr
\noalign{\vskip-10pt}
}
}
}
\endinsert

These two lists are Tables~I and~II, essentially.  I have taken
the liberty of modifying Berger's lists slightly, dropping the entries
in the original list that did not actually satisfy Berger's two criteria 
($\SO^*(2n)\simeq\SO(n,\bbH)$, which does not satisfy the first criterion%
\footnote{$^2$}{\ninerm An observation due to R. McLean.}%
, and the $\Spin(9,\bbC)$-type entries, which do not satisfy the second 
criterion%
\footnote{$^3$}{\ninerm Observed independently by Alekseevskii~[Al1] and 
Brown and Gray~[BG].}%
) and including $\Sp(p,\bbR){\cdot}\SL(2,\bbR)$, the inadvertently 
omitted `split form' of the quaternionic K\"ahler case.  
In the nonmetric case, I have amplified the list somewhat by making 
explicit the various real forms as well as the fact that, except 
for~$\CO(p,q)$, each of the entries on Berger's nonmetric list represents 
only the semi-simple part of~$H$, to which one can add an arbitrary subgroup 
of the (abelian) commuting subgroup to make the full group~$H$.

{\it Exotic Holonomies.}
The nonmetric list supplied by Berger was never claimed to be
a complete list, though it was supposed to have omitted at most a
finite number of possibilities.  
The full list of omissions, comprising Tables~III and~IV and nowadays 
referred to as the {\it exotic list\/}, was recently compiled by a 
combination of the efforts of Chi, Merkulov, Schwachh\"ofer, and myself.  
This will be reported on in \S3, along with the reasons for the
division into two lists.

\topinsert
\centerline{
\vtop{\offinterlineskip
\halign{
&\vrule#&\strut\quad\hfil#\hfil\quad\cr  %The preamble
\multispan4{\hfil\bf III. Exotic Conformal Holonomies\vphantom{p}\hfil}\cr
\noalign{\smallskip}
\multispan4{\hrulefill}&\cr
height 2 pt&\omit&&\omit&\cr
& H$^\lozenge$ && $\eum$ &\cr
height 2 pt&\omit&&\omit&\cr
\multispan4{\hrulefill}&\cr
height 1 pt&\omit&&\omit&\cr
\multispan4{\hrulefill}&\cr
height 2 pt&\omit&&\omit&\cr
& $\bbR^+\!\cdot\!\SL(2,\bbR)$ &
    & $\bbR^4\simeq\Sym^3(\bbR^2)$ &\cr
height 1 pt&\omit&&\omit&\cr
& $\bbC^*\!\cdot\!\SL(2,\bbC)$ &
    & $\bbC^4\simeq\Sym^3(\bbC^2)$ &\cr
height 2 pt&\omit&&\omit&\cr
\multispan4{\hrulefill}&\cr
height 2 pt&\omit&&\omit&\cr
& $G_\bbC\!\cdot\!\SL(2,\bbR)\,^a$ &
    & $\bbC^2\simeq\bbC{\otimes_\bbR}\bbR^2$ &\cr
height 1 pt&\omit&&\omit&\cr
& $G_\bbC\!\cdot\!\Sp(1)\,^b$      &
    & $\bbC^2\simeq\bbC{\otimes_\bbC}\bbH\phantom{^2}$ &\cr
height 2 pt&\omit&&\omit&\cr
\multispan4{\hrulefill}&\cr
height 2 pt&\omit&&\omit&\cr
& $G_\bbR\!\cdot\!\Spin(5,5)$ && $\bbR^{16}$ &\cr
height 1 pt&\omit&&\omit&\cr
& $G_\bbR\!\cdot\!\Spin(1,9)$ && $\bbR^{16}$ &\cr
height 1 pt&\omit&&\omit&\cr
& $G_\bbC\!\cdot\!\Spin(10,\bbC)$ && $\bbC^{16}$ &\cr
height 2 pt&\omit&&\omit&\cr
\multispan4{\hrulefill}&\cr
height 2 pt&\omit&&\omit&\cr
& $G_\bbR\!\cdot\!\E_6^1$ && $\bbR^{27}$ &\cr
height 1 pt&\omit&&\omit&\cr
& $G_\bbR\!\cdot\!\E_6^4$ && $\bbR^{27}$ &\cr
height 1 pt&\omit&&\omit&\cr
& $G_\bbC\!\cdot\!\E_6^\bbC$ && $\bbC^{27}$ &\cr
height 2 pt&\omit&&\omit&\cr
\multispan4{\hrulefill}&\cr
\noalign{\smallskip}
\omit\ $^\lozenge$ {\ninerm $G_\bbF$ is any connected subgroup of~$\bbF^*$}
\hidewidth\cr
\noalign{\vskip5pt}
\omit Restrictions: \hidewidth\cr
\noalign{\vskip7pt}
\omit\ \ \ $^a$ {\ninerm $G_\bbC\not\subseteq\bbR^*$ (for irreducibility).}
\hidewidth\cr
\noalign{\vskip4pt}
\omit\ \ \ $^b$ {\ninerm $G_\bbC\not\subseteq S^1$ (to be nonmetric).}
\hidewidth\cr
\noalign{\vskip4pt}
\omit\ \ \ $^c$ {\ninerm $p+q\ge3$ (for irreducibility).}\hidewidth\cr
\noalign{\vskip4pt}
\omit\ \ \ $^d$ {\ninerm $n\ge3$ (for irreducibility).}\hidewidth\cr
\noalign{\vskip4pt}
\omit\ \ \ $^e$ {\ninerm $n\ge2$ (to be nonmetric).}\hidewidth\cr
}
\vskip-10pt}
\hskip1pt
\vtop{\offinterlineskip
\halign{
&\vrule#&\strut\quad\hfil#\hfil\quad\cr  %The preamble
\multispan4{\hfil \bf IV. Exotic Symplectic Holonomies \hfil}\cr
\noalign{\smallskip}
\multispan4{\hrulefill}&\cr
height 2 pt&\omit&&\omit&\cr
& H && $\eum$ &\cr
height 2 pt&\omit&&\omit&\cr
\multispan4{\hrulefill}&\cr
height 1 pt&\omit&&\omit&\cr
\multispan4{\hrulefill}&\cr
height 2 pt&\omit&&\omit&\cr
& $\SL(2,\bbR)$ && $\bbR^4\simeq\Sym^3(\bbR^2)$ &\cr
height 1 pt&\omit&&\omit&\cr
& $\SL(2,\bbC)$ && $\bbC^4\simeq\Sym^3(\bbC^2)$ &\cr
height 2 pt&\omit&&\omit&\cr
\multispan4{\hrulefill}&\cr
height 2 pt&\omit&&\omit&\cr
& $\SL(2,\bbR)\!\cdot\!\SO(p,q)\,^c$ && $\bbR^{2}\otimes\bbR^{p+q}$ &\cr
height 1 pt&\omit&&\omit&\cr
& $\SL(2,\bbC)\!\cdot\!\SO(n,\bbC)\,^d$ && $\bbC^{2}\otimes\bbC^{n}$ &\cr
height 1 pt&\omit&&\omit&\cr
& $\Sp(1)\!\cdot\!\SO(n,\bbH)\,^e$ && $\bbH^{n}$ &\cr
height 2 pt&\omit&&\omit&\cr
\multispan4{\hrulefill}&\cr
height 2 pt&\omit&&\omit&\cr
& $\Sp(3,\bbR)$ && $\bbR^{14}\simeq \Lam^3_0(\bbR^6)$ &\cr
height 1 pt&\omit&&\omit&\cr
& $\Sp(3,\bbC)$ && $\bbC^{14}\simeq \Lam^3_0(\bbC^6)$ &\cr
height 2 pt&\omit&&\omit&\cr
\multispan4{\hrulefill}&\cr
height 2 pt&\omit&&\omit&\cr
& $\SL(6,\bbR)$ && $\bbR^{20}\simeq\Lam^3(\bbR^6)\phantom{^\bbR}$ &\cr
height 1 pt&\omit&&\omit&\cr
& $\SU(1,5)$ && $\bbR^{20}\simeq\Lam^3(\bbC^6)^\bbR$ &\cr
height 1 pt&\omit&&\omit&\cr
& $\SU(3,3)$ && $\bbR^{20}\simeq\Lam^3(\bbC^6)^\bbR$ &\cr
height 1 pt&\omit&&\omit&\cr
& $\SL(6,\bbC)$ && $\bbC^{20}\simeq\Lam^3(\bbC^6)\phantom{^\bbR}$ &\cr
height 2 pt&\omit&&\omit&\cr
\multispan4{\hrulefill}&\cr
height 2 pt&\omit&&\omit&\cr
& $\Spin(2,10)$ && $\bbR^{32}$ &\cr
height 1 pt&\omit&&\omit&\cr
& $\Spin(6,6)$ && $\bbR^{32}$ &\cr
height 1 pt&\omit&&\omit&\cr
& $\Spin(6,\bbH)$ && $\bbR^{32}$ &\cr
height 1 pt&\omit&&\omit&\cr
& $\Spin(12,\bbC)$ && $\bbC^{32}$ &\cr
height 2 pt&\omit&&\omit&\cr
\multispan4{\hrulefill}&\cr
height 2 pt&\omit&&\omit&\cr
& $\E_7^5$ && $\bbR^{56}$ &\cr
height 1 pt&\omit&&\omit&\cr
& $\E_7^7$ && $\bbR^{56}$ &\cr
height 1 pt&\omit&&\omit&\cr
& $\E_7^\bbC$ && $\bbC^{56}$ &\cr
height 2 pt&\omit&&\omit&\cr
\multispan4{\hrulefill}&\cr
}
\vskip-10pt}
}
\endinsert

\vskip 0.5truecm
\noindent
{\bf 1.6.} {\smc Local Existence.}
Berger's lists (suitably modified) provide possibilities for
irreducibly acting holonomy groups, but to verify that these possibilities 
actually can occur requires methods beyond representation theory.  Most of 
the methods that have been employed can be grouped into a small number of 
categories:
\smallskip
{\it Explicit construction.}  This is the simplest method, when it is
available.  For the metric list, there are locally symmetric examples
with every holonomy except the special K\"ahler cases, where the 
holonomy is~$\SU(p,q)\subset\GL(\bbC^{p+q})$; the hyper-K\"ahler cases, 
where the holonomy is $\Sp(p,q)\subset\GL(\bbH^{p+q})$; and the 
`exceptional' holonomies, which comprise the 
groups~$\G_2^\bbC\subset\GL(\bbC^7)$ and~$\Spin(7,\bbC)\subset\GL(\bbC^8)$ 
and certain of their real forms.  Of course, one would like to know that
the locally symmetric examples are not the only ones.

Sometimes constructing examples is easy:  The generic 
pseudo-Riemannian metric has holonomy is~$\SO(p,q)$.

In other cases, simple underlying geometric structures can be used as a 
starting point.  For example,
all complex structures are flat, and the general (pseudo-)K\"ahler metric
can be described in the standard background complex structure by means 
of a (pseudo-)K\"ahler potential. For the generic choice of such a
potential, the holonomy will be~$\U(p,q)\subset\SO(2p,2q)$.
Another example is the special K\"ahler case, where one can start with a 
background complex structure with a specified holomorphic volume form 
and then find a K\"ahler metric preserving this volume form by requiring 
that the K\"ahler potential satisfy a single second order, elliptic equation.

One can also find examples by looking for those with a large symmetry group.
For the hyper-K\"ahler case, one can start with the complex symplectic 
structure on the cotangent bundle of certain Hermitian symmetric spaces 
and look for a K\"ahler potential compatible with this complex symplectic 
structure that is invariant under the action of the isometry group,
thereby reducing the problem to solving an ordinary differential equation.
This was Calabi's method for constructing a hyper-K\"ahler metric 
on~$T^*\bbC\bbP^n$, the first known example in general dimensions.
In the quaternionic K\"ahler case, where the holonomy 
is~$\Sp(p,q){\cdot}\Sp(1)\subset\GL(\bbR^{4(p+q)})$, Alekseevskii found
homogeneous nonsymmetric examples on certain solvable Lie groups.
Even for the exceptional holonomies, there are
explicit examples of cohomogeneity one~[Br1], [BS].

Examples in the hyper-K\"ahler and quaternionic K\"ahler cases can also
be constructed by the method of {\it reduction\/}, which takes advantage of 
descriptions of these structures in terms of multi-symplectic geometry,
generalizing the well-known method of Marsden-Weinstein 
reduction in symplectic geometry so as to handle the multi-symplectic case.
For an account, see~[Bes, Addendum~E].

\smallskip
{\it Twistor Methods.} After Penrose's description of the self-dual
metrics in dimension~$4$, Hitchin and Salamon~[Sa], among others, 
were able to generalize this method to describe the hyper-K\"ahler and 
quaternionic K\"ahler metrics in terms of natural holomorphic geometric 
structures on the moduli space of rational curves with certain simple
normal bundles in a complex manifold.

In fact, it was the study~[Br2] of the moduli space of rational
curves on a complex surface with normal bundle~${\cal O}(3)$ that 
turned up the first known examples (the first two entries in each of 
Tables~III and~IV) of omissions%
\footnote{$^4$}{\baselineskip=9pt{\ninerm 
These examples were referred to as `exotic' in~[Br2] and the term has been 
adopted to describe any nonmetric subgroup~$H\subset\GL(\eum)$ that satisfies 
Berger's criteria but that does not appear on Berger's original 
nonmetric list.}\par}
from Berger's nonmetric list.
Moreover, in the holomorphic category, it was shown that any connection 
with one of these holonomies could be constructed as the natural connection 
on the four-dimensional component of the moduli space of Legendrian rational 
curves in a complex contact three-fold.

Following this, Merkulov~[Me1] showed that this approach could be 
generalized to cover the geometry of the moduli space of Legendrian 
deformations of certain complex homogeneous spaces and began to discover
more exotic examples.  This and its further developments will be reported 
on in \S3.

\smallskip
{\it Exterior differential systems.}  Another approach is to treat
the equation~$\bar T=0$ directly as a system of PDE
for sections of the bundle~$F(M,\eum)/H\to M$.  In nearly all cases, 
this method leads to the study of an overdetermined system of PDE,
so that Cartan-K\"ahler machinery must be brought to bear.  For
definitions and results regarding Cartan-K\"ahler theory, the reader 
can consult~[BCG].
  
The general approach can be summarized as follows:  Consider the 
structure equations derived so far for a torsion-free connection~$\theta$ 
on an~$H$-structure~$B\subset F(M,\eum)$,
$$
\eqalign{
d\omega &= -\theta\w\omega\cr
d\theta &= -{\ts{1\over2}}[\theta,\theta] 
             + {\ts{1\over2}}\,R(\omega\w\omega)\cr
dR &= -\theta.R + R'(\omega)\cr
}
\leqno{(0.4)}
$$
where~$R: B\to K(\euh)$ and ~$R':B\to K^1(\euh)\subset K(\euh)\otimes\eum^*$
are as defined before.  There are two things that need to be checked in order
to be able to apply Cartan's general existence theorem for coframings
at this level:  First, the inclusion~$K^1(\euh)\subset K(\euh)\otimes\eum^*$
should be an involutive tableau in Cartan's sense.  Second, there should be
a quadratic map~$Q:K(\euh)\to K^1(\euh)\otimes\eum^*$ so that the exterior
derivative of the third structure equation in (0.4) can be written
in the form~$\bigl(dR' + \theta.R' - Q(R)(\omega)\bigr)(\omega) = 0$.
(This is the familiar `vanishing torsion' condition in exterior differential
systems.)  

When these two conditions are satisfied, Cartan's existence
theorem asserts that, up to local diffeomorphism, the real analytic 
torsion-free connections with holonomy lying in~$H\subset\GL(\eum)$ depend
on~$s_q$ functions of~$q$ variables, where~$s_q$ is the last nonzero
Cartan character of~$K^1(\euh)\subset K(\euh)\otimes\eum^*$.  Moreover,
for any~$(R_0,R'_0)\in K(\euh)\times K^1(\euh)$, there exists a 
torsion-free connection~$\theta$ on an~$H$-structure~$B\subset F(\eum,\eum)$ 
and a~$u_0\in B$ for which the curvature functions~$R$ and $R'$ satisfy
$R(u_0)= R_0$ and $R'(u_0)=R'_0$. 

When one can choose the 
element~$R_0\in K(\euh)\subset\euh\otimes\Lam^2(\eum^*)$ so that it is
surjective%
\footnote{$^5$}{\baselineskip=9pt{\ninerm Actually, one only needs that 
the image~$R_0\bigl(\Lam^2(\eum)\bigr)\subset\euh$ generates~$\euh$ as 
an algebra.}\par}
 as a map~$R_0:\Lam^2(\eum)\to \euh$, it will then follow from the
Ambrose-Singer holonomy theorem that the holonomy of~$\theta$ 
contains the identity component of~$H$.  If, moreover, one can choose
$R'_0$ to be nonzero, such a connection~$\theta$ will not be locally
symmetric.

This analysis applies successfully to each of the entries of 
Tables~I,~II, and~III.  By contrast, for each of the
entries of Table~IV, the tableau~$K^1(\euh)\subset K(\euh)\otimes\eum^*$
is not involutive.  Further discussion of this point will
be given in~\S3.

Generally, this method is good only for local analysis, but it has the 
distinct advantage that it not only proves existence of connections
with a given holonomy, but provides their `degree of generality' in Cartan's 
sense.  For example, Table~I gives the degree of generality of each of the 
possible pseudo-Riemannian, irreducible holonomies.  For a similar 
discussion of the nonmetric list, see the survey~[Br3], where various
simplifications of the general argument are introduced to shorten the
exposition.

This method was first used to prove the existence of metrics with 
holonomy~$\G_2$ and $\Spin(7)$ and is still the only method that 
constructs the general local solution and describes its 
degree of generality.  This is also the only known method for analyzing
Entries~3 and~4 of Table~III.

\smallskip
{\it Poisson Constructions.} The examples~$H$ in Table IV are subgroups
of~$\Sp(S)\subset\GL(\eum)$ for a nondegenerate skewsymmetric bilinear
form~$S$ on~$\eum$.  Hence the corresponding $H$-structures (when they
exist) have an underlying symplectic structure.

For the first two examples from Table~IV, each torsion-free
connection~$(M,B,\theta)$ of these types was analyzed and reconstructed
in~[Br2] from its derived curvature 
map~$J=(R,R'):B\to K(\euh)\oplus K^1(\euh)$.   This reconstruction involved
a number of seemingly miraculous identities, but since only
these two examples were known, it did not seem worthwhile to look for 
an interpretation of these identities. However, when Chi, Merkulov, and
Schwachh\"ofer~[CMS] found other exotic symplectic examples, they noticed
that this reconstruction technique generalized and they were able to
explain it in the context of Poisson geometry in a very beautiful
way.  This, too, will be reported on in~\S3.

\vskip 0.5truecm
\noindent
{\bf 1.7.} {\smc Compact Riemannian Examples.}  The history of compact 
Riemannian manifolds with reduced holonomy groups is long and complex, 
so I will not attempt a full account here.  For more details on the cases 
I only mention, the reader can consult the relevant chapters of~[Bes]
and the references cited therein.
\smallskip
{\it K\"ahler manifolds.}  This subject has the longest history, 
predating even Berger's classification.  Every smooth algebraic variety
carries a K\"ahler structure and this accounts for their importance in
algebraic geometry.
\smallskip
{\it Special K\"ahler manifolds.}  The major milestone here is, 
of course, Yau's solution in the mid 1970s of the Calabi conjecture, showing 
that every compact K\"ahler manifold with trivial canonical bundle carries 
a special K\"ahler structure.  For this reason, compact manifolds 
endowed with such a structure are usually referred to as Calabi-Yau
manifolds.
\smallskip
{\it Hyper-K\"ahler manifolds.}  In the early 1980s, Beauville
and Mukai were each able to construct compact, simply connected
K\"ahler manifolds~$M^{4p}$ that carried a nondegenerate complex 
symplectic form (Fujiki had constructed examples for~$p=2$ slightly
earlier).  These necessarily had trivial canonical bundle, so by Yau's 
solution of the Calabi conjecture, these carried special K\"ahler structures.
By an argument of Bochner, the complex symplectic form had to be parallel
with respect to this special K\"ahler structure, which forced the
holonomy to lie in~$\Sp(p)$.  Further arguments showed that these 
examples were not products of lower dimensional complex manifolds
and this implied that the holonomy had to actually be~$\Sp(p)$.
These were the first known compact examples.  

\smallskip
{\it  Quaternion K\"ahler manifolds.}  In this case, all known 
compact examples are locally symmetric, but we know of no
reason why this should be true, except for dimension~$8$ [PoS].  
Of course, a great deal is known about the possible geometry of 
such examples, see~[Sa].

\smallskip
{\it  $\G_2$ and $\Spin(7)$ manifolds.}  The remarkable recent
work of Joyce~[Jo1,2] establishes the existence of compact manifolds
with these holonomies.  This will be reported on in \S2.

\vskip 1.2truecm
\noindent
{\bf 2. R\ninebf IEMANNIAN HOLONOMY}

In this section, $g$ will denote a smooth Riemannian metric on a 
connected, smooth manifold~$M^n$.  The reference space~$\eum$ will
be taken to be~$\bbR^n$ with its standard inner product, and 
$\O(\eum)\subset\GL(\eum)$ will denote its orthogonal group. 
The~$\O(\eum)$-structure~$B\subset F(M,\eum)$ consisting of the 
coframes~$u:T_xM\to\eum$ that are isometries of vector spaces is the 
orthonormal coframe bundle and the Levi-Civita connection on~$B$ will 
be denoted~$\theta$.  This is, of course, the unique torsion-free 
$\O(\eum)$-connection on~$B$.  When there is no danger of
confusion, I will simply write~$H_u$ and $B_u$ instead of~$H^\theta_u$
and~$B^\theta_u$.

One thing that makes the Riemannian case simpler to deal with 
than others is the {\it de~Rham Splitting Theorem\/}~[Bes, KNo], 
which occurs in a local form and a global form.  
The local form asserts that if~$(H_u)^0$ acts reducibly on~$\eum$,
say, preserving irreducible orthogonal 
subspaces~$\eum_i\subset\eum$ for~$1\le i\le k$,
then $(H_u)^0$ is the direct product of its 
subgroups~$(H_u)_i^0$, where~$(H_u)_i^0$ is the subgroup that acts trivially
on~$\eum_j$ for~$j\not=i$.  Moreover, the metric~$g$ locally splits as a
product in a corresponding way%
\footnote{$^6$}{\baselineskip=9pt{\ninerm
What is also true, but not obvious, is that each of the groups
$(H_u)_i^0\subset\SO(\eum_i)$ is the holonomy of some Riemannian metric,
even when it is not the holonomy of the corresponding local factor 
of~$g$. See [Bes, Theorem 10.108]}\par}.  
The global form asserts that if, in
addition, $M$ is simply connected and the metric~$g$ is complete, then
$M$ globally splits as a Riemannian product
$$
(M,g) = (M_1,g_1)\times\cdots\times(M_k,g_k).
$$
so that~$(H_u)^0_i = (H_u)_i$ is the holonomy of~$(M_i,g_i)$.

\vskip 0.5truecm
\noindent
{\bf 2.1.} {\smc Non-closed holonomy.}  While the remainder of this
report deals only with connected holonomy groups, I cannot pass up the
opportunity to mention a recent result of particular interest in
Riemannian holonomy.  It had been a question for some time whether the
holonomy of a compact Riemannian manifold must necessarily be 
compact.  Using the de Rham Splitting Theorem and Berger's holonomy 
classification in the irreducible Riemannian case, one sees that this
is so if~$M$ is simply connected.  Thus, for any Riemannian manifold,
the restricted holonomy group is compact, so it becomes a question of
whether the fundamental group can cause the holonomy group to have
an infinite number of components even when~$M$ is compact. 

Very recently, B. Wilking~[Wi] has shown that this can indeed occur.
He has produced an example of a compact manifold whose holonomy group 
does have an infinite number of components.  His example is of the form~$M^5 
=\Gamma\backslash(\bbR^2{\times}N^3)$ where $\Gamma$ is a subgroup of the 
isometry group of~$\bbR^2{\times}N^3$ that acts properly discontinuously
and cocompactly.

\vskip 0.5truecm
\noindent
{\bf 2.2.}  {\smc Compact manifolds with exceptional holonomy.}  The most
remarkable development in Riemannian holonomy in recent years has been
the spectacular construction by Dominic Joyce of compact $7$-manifolds
with holonomy~$G_2$ and compact $8$-manifolds with holonomy $\Spin(7)$.
His constructions are full of new ideas and, while it is not difficult
to outline these ideas, their successful execution turns out to require
very careful, subtle estimates.  I will not attempt to explain these, but
refer the reader to the original sources~[Jo1,2] and to the forthcoming
book~[Jo3].  Also, I will concentrate on the $\G_2$ case, as the 
$\Spin(7)$ case follows the same spirit, but the details are different.
\smallskip
{\it The fundamental $3$-form.}  Let~$\eum=\bbR^7$.  It has been
known for some time~[Br1] that there is an open 
$\GL(\eum)$-orbit~$\Lam_+^3(\eum^*)$ in the 3-forms on~$\eum$ 
so that the stabilizer of any element~$\phi\in\Lam_+^3(\eum^*)$ is 
a compact connected simple Lie group of dimension~$14$ and which is
therefore isomorphic to~$\G_2$.  Consequently,
for any $7$-manifold~$M$, there is an open subbundle~$\Lam_+^3(T^*M)
\subset\Lam^3(T^*M)$ so that the $\G_2$-structures on~$M$ are in
one-to-one correspondence with the sections~$\Omega_+^3(TM)$ of this bundle.  
Such a section will exist if and only if $M$ is orientable and
spinnable~[LM], so I assume this from now on.

Associated to any section~$\sigma\in\Omega_+^3(M)$, there is a
canonical metric~$g_\sigma$ and orientation, inducing a well-defined
Hodge star operator~$\ast_\sigma:\Omega^p(M)\to\Omega^{7-p}(M)$.  If 
$\sigma$ is parallel with respect to the Levi-Civita connection of~$g_\sigma$
then the holonomy of~$g_\sigma$ will be a subgroup of~$\G_2$ and the
$\G_2$-structure associated to~$\sigma$ will be torsion-free.  By
a result of Gray, this~$\G_2$-structure is torsion-free
if and only if both~$\sigma$ and~$\ast_\sigma\sigma$ are closed.
Conversely, if a Riemannian metric~$g$ on~$M^7$ has holonomy a subgroup of
$\G_2$, then there will be a $g$-parallel section~$\sigma\in\Omega_+^3(M)$
for which~$g=g_\sigma$.

A metric with holonomy a subgroup of~$\G_2$ is known to be
Ricci flat~[Bes, 10.64], so if~$M$ is compact, Bochner vanishing shows
that its harmonic 1-forms and its Killing fields must be parallel.  
Thus, in the compact case,
$b_1(M)=0$ if the holonomy actually equals all of~$\G_2$.  Combining
this information with the Cheeger-Gromoll splitting theorem~[Bes, 6.67],
one finds that a compact manifold with holonomy~$\G_2$ must have finite
fundamental group, so one can, without loss of generality, assume that
$M$ is simply connected, which I do from now on.  Conversely, a
compact simply connected Riemannian $7$-manifold whose holonomy is a 
subgroup of~$\G_2$ must actually equal $\G_2$ since any proper subgroup 
of~$\G_2$ that can be a holonomy group fixes a nontrivial subspace.

One can now consider the moduli space~$\cal M$ consisting of the
sections~$\sigma\in\Omega_+^3(M)$ satisfying~$d\sigma=d\,{\ast_\sigma}\sigma=0$
modulo diffeomorphisms of~$M$ isotopic to the identity.  There is
an obvious `Torelli' map~$\tau: {\cal M}\to H^3_{\rm dR}(M)$ defined by
$\tau\bigl([\sigma]\bigr) = [\sigma]_{\rm dR}$.  Joyce's first striking
result~[Jo1] is the analog of the local Torelli theorem, namely that $\tau$ 
is locally one-to-one and onto, in fact, a local diffeomorphism in the
natural smooth structure on~$\cal M$.  Thus, the moduli space is 
said to be `unobstructed'.  

Next, Joyce proves a remarkable existence theorem:  
If~$\sigma\in\Omega_+^3(M)$ is closed, then there is a 
constant~$C$ that depends on the norm of the curvature of the
metric~$g_\sigma$, its volume, and its injectivity radius, so that,
if~$\bigl|d(\ast_\sigma\sigma)\bigr|_\sigma<C$, then there exists an 
exact 3-form~$\phi$ so that $\sigma+\phi$ lies 
in~$\Omega_+^3(M)$ and is closed and coclosed.  In other words, a closed 
3-form in~$\Omega_+^3(M)$ that is `close enough' to being coclosed can 
be perturbed to a cohomologous 3-form in~$\Omega_+^3(M)$ that is both 
closed and coclosed.

Thus, to prove the existence of a compact Riemannian $7$-manifold with 
holonomy~$\G_2$, it suffices to construct a simply connected $7$-manifold
endowed with a closed $3$-form that satisfies such a `close enough' estimate.

Joyce's idea for doing this is extremely clever:  He starts with the
flat~$\G_2$-structure~$\sigma_0$ on the $7$-torus~$T^7 = \bbR^7/\bbZ^7$
and passes to a simply connected quotient orbifold~$X=\Gamma\backslash T$ 
where~$\Gamma$ is a finite group of~$\sigma_0$-symmetries.  This provides a 
flat $\G_2$-orbifold whose singular locus is 
a finite number of 3-tori~$T^3$, each of which has a neighborhood 
of the form~$T^3\times B^4/\{\pm\}$
where~$B^4/\{\pm\}$ is the standard $4$-ball around the origin 
in~$\bbR^4=\bbC^2$ divided by the equivalence relation~$v\sim-v$.  

Now, it has been known
for a long time that $\bbR^4/\{\pm\}$ is metrically the scaling limit of
the~$\SU(2)$-holonomy metric on~$T^*\bbC\bbP^1$ constructed by
Eguchi and Hanson as one scales the metric to contract the zero section
to a point.  Because~$I_3{\times}\SU(2)\subset\SO(7)$ is a subgroup 
of~$\G_2$, it follows that one can regard the flat $\G_2$-structure
on~$X$ in each singular locus neighborhood~$T^3\times B^4/\{\pm\}$ as
the limit of a scaling of a $\G_2$-structure on~$T^3{\times}T^*\bbC\bbP^1$.
Joyce cuts out these singular neighborhoods and glues back in $T^3\times N$
where~$N$ is a neighborhood of the zero section in~$T^*\bbC\bbP^1$, smoothly
joining the flat $3$-form with the Eguchi-Hanson-derived (closed)
$3$-form on the overlaps. This `surgery' produces a smooth manifold~$\hat X$,
but does not disturb the fundamental group, which remains trivial. 

By being very careful (here is where Joyce's estimates
are extremely delicate), he shows that he can do this in such a way 
that the resulting closed 3-form~$\sigma\in\Omega_+^3(\hat X)$ 
satisfies his estimate.  I.e., it is close enough to being coclosed 
that it can be perturbed to a $\tilde\sigma\in\Omega_+^3(\hat X)$ that 
is both closed and coclosed.  Of course, since~$X$ is simply connected, 
it follows that the resulting torsion-free $\G_2$-structure has holonomy 
equal to~$\G_2$.

By applying this idea to a number of different finite groups~$\Gamma$, 
Joyce has been able to construct $\G_2$-metrics on a number of different
$7$-manifolds.

A similar set of ideas allows Joyce to construct compact 8-manifolds
with holo\-nomy~$\Spin(7)$, once one realizes that~$\Spin(7)$ can
be defined in~$\GL(8,\bbR)$ as the stabilizer of a certain $4$-form
in eight variables.  The interested reader should consult [Jo1,2].

\vskip 1.2truecm
\noindent
{\bf 3. I\ninebf RREDUCIBLE TORSION-FREE NON-METRIC CONNECTIONS}

\vskip 0.5truecm
\noindent
{\bf 3.1.} {\smc Twistor constructions.}  As previously mentioned,
I found the first examples of `exotic' holonomy groups by studying
the geometry of the moduli space~$M$ of rational curves (i.e., complex
curves of genus~$0$) in a complex surface~$S$ with normal 
bundle~${\cal O}(3)$.  Following the examples of Hitchin's study of
rational curves in a surface~$S$ with normal bundle~${\cal O}(k)$ 
for~$k=1$ and~$2$, I knew that~$M$ would have dimension~$4$
and would have a natural~$G_3$-structure, where~$G_3\subset\GL(4,\bbC)$
is the image of~$\GL(2,\bbC)$ acting by linear substitutions on
the $4$-dimensional space~$V_3$ of cubic polynomials in two variables.
I also knew that there would be a canonical $G_3$-connection from 
general principles, but I was very surprised to find that this 
connection was torsion-free.  

In the examples Hitchin had analyzed, the geometry on the moduli
space allowed one to reconstruct the surface~$S$ and so I fully 
expected to be able to do the same in this case.  However, it turned 
out that the story was more subtle than that.  In the standard 
double fibration picture:
$$
\matrix{
&&I\hskip-.5em&&\cr
&\raise.5em\hbox{$\lambda$}\hskip-.5em\swarrow\hskip-.5em&
&\searrow\hskip-.5em\raise.5em\hbox{$\rho$}\hskip-1em&\cr
M\hskip-1em&&&&S\cr }
$$
where~$I\subset M\times S$ is the set of pairs~$(C,p)$ where $p\in S$
is a point of the rational curve~$C\in M$ and~$\lambda$ and~$\rho$
are just the projections onto the factors, each~$p\in S$ would
correspond to a hypersurface~$H_p=\lambda\bigl(\rho^{-1}(p)\bigr)$ in~$M$
(since it is only one condition for a curve in~$S$ to pass through
a given point~$p$).  The members of this $2$-parameter family of 
hypersurfaces in~$M$ would be expected to be the solutions of some
differential geometric problem in~$M$, but I was not able to find
a geometrically defined $2$-parameter family of hypersurfaces in
the general $4$-manifold carrying a torsion-free $G_3$-structure.  

However, a $3$-parameter family~$Y$ of $2$-dimensional surfaces did
present itself.  This can be described as follows:  By the defining
property of a~$G_3$-structure~$B$ on~$M^4$, each tangent space~$T_xM$ can be 
thought of as the space of homogeneous cubic polynomials in two variables.
This defines a distinguished~$\bbP^1$ of lines, namely the perfect 
cubes, and a distinguished~$\bbP^1$ of $2$-planes, namely the multiples
of a perfect square.  I called such lines and 2-planes {\it null\/}.
It was not difficult to show that, when~$B$ had a torsion-free connection,
each null $2$-plane was tangent to a unique totally geodesic $2$-surface 
in~$M$ all of whose tangent planes were null.  This family~$Y$ then
fit into a double fibration%
\footnote{$^7$}{\baselineskip=9pt{\ninerm 
Assuming~$Y$ to be Hausdorff in its natural topology, 
which can always be arranged by replacing~$M$ by a $\theta$-convex 
open set in~$M$.}\par}
$$
\matrix{
&& N^5\hskip-.5em&&\cr
&\raise.5em\hbox{$\lambda$}\hskip-.5em\swarrow\hskip-.5em&
&\searrow\hskip-.5em\raise.5em\hbox{$\rho$}\hskip-1em&\cr
M^4\hskip-1em&&&&Y^3\cr }
$$
where the fibers of~$\lambda$ are~$\bbP^1$'s.  Moreover, I was able
to show that~$Y$ carried a natural structure as a contact manifold,
that the family of~$\bbP^1$'s given by~$C_x=\rho\bigl(\lambda^{-1}(x)\bigr)$ 
for~$x\in M$ were all Legendrian curves for this contact structure,
and that, moreover, this family was an open set in the space of
Legendrian rational curves in~$Y$.  (In the surface case that I had
started out with, $Y$ turned out to be the projectivized tangent bundle
of~$S$.)  

I then showed that the picture could be reversed:  Starting with
a holomorphic contact $3$-manifold~$Y$, one could look at the 
moduli space~$M$ of rational Legendrian curves~$C$ in~$Y$ to which
the contact bundle~$L\subset T^*Y$ restricts to be isomorphic to
${\cal O}(-3)$ and show that it was a smooth moduli space of 
dimension~$4$ on which there was a canonical torsion-free $G_3$-structure.

{\it Merkulov's generalization.} 
In a remarkable series of papers, Merkulov~[Me1,2,3] showed that this
moduli space and double fibration construction 
obtains in a very general setting, starting from the data of an irreducibly 
acting (and therefore reductive) complex subgroup~$H\subset\GL(n,\bbC)$, 
a complex $n$-manifold~$M$, and a holomorphic\break
$H$-structure~$B\subset F(M,\bbC^n)$ endowed with a holomorphic
torsion-free connection~$\theta$.

The semi-simple part~$H_s\subset H$ acts irreducibly on~$\bbC^n$.
With respect to a Cartan subalgebra of~$H_s$ and an ordering of its roots, 
there will be a unique line~$E\subset(\bbC^n)^*$ spanned by a vector of 
highest weight.  The $H_s$-orbit~$F\subset\bbP\bigl((\bbC^{n})^*\bigr)$ 
of~$E$ is a minimal $H_s$-orbit in $\bbP\bigl((\bbC^{n})^*\bigr)$, 
a generalized flag variety of~$H_s$ of some dimension~$k\le n{-}1$,
endowed with the hyperplane section bundle~${\cal L}$.
(In the original case I treated, $F$ is the projectivization 
of the set of perfect cubes and hence is a $\bbP^1$.  The bundle~${\cal L}$
is ${\cal O}(-3)$.)  
The $H$-structure~$B$ provides identifications
$T_xM\simeq\bbC^n$ unique up to an action of~$H$,
so there is a subbundle~$N\subset\bbP(T^*M)$
whose fiber~$N_x\subset\bbP(T_x^*M)$ over~$x$ corresponds 
to~$F\subset\bbP\bigl((\bbC^{n})^*\bigr)$ via any $B$-identification.  

The projectivized cotangent bundle of any manifold is 
canonically a contact manifold, and the torsion-free
condition on~$\theta$ immediately implies that~$N$ is an 
{\it involutive\/} submanifold of~$\bbP(T^*M)$, i.e.,
the restriction of the contact structure to~$N$ is degenerate, with
Cauchy leaves of the largest possible dimension, namely $n{-}k{-}1$,
the codimension of~$N$ in~$\bbP(T^*M)$.
When it is Hausdorff
\footnote{$^8$}{\baselineskip=9pt{\ninerm 
This can always be arranged by replacing~$M$ by a suitably $\theta$-convex 
open set in~$M$.}\par},
 the leaf space~$Y$ of this Cauchy foliation is 
canonically a contact manifold, yielding the double fibration
$$
\matrix{
&& N^{n+k}\hskip-.5em&&\cr
&\raise.5em\hbox{$\lambda$}\hskip-.5em\swarrow\hskip-.5em&
&\searrow\hskip-.5em\raise.5em\hbox{$\rho$}\hskip-1em&\cr
M^n\hskip-1em&&&&Y^{2k+1}\cr }
$$
where the manifolds~$F_x 
= \rho\bigl(\lambda^{-1}(x)\bigr)$ are Legendrian $k$-dimensional
submanifolds of~$Y$.  

Merkulov then goes on to prove that, nearly always, one can
recover~$M$ as the complete moduli space of the Legendrian 
immersions~$F\subset Y$ that pull back the contact bundle~$L\subset T^*Y$ 
to be~$\cal L$. Moreover, when~$H_s$ acts as 
the full biholomorphism group of~$F$ (which, again, is nearly always) 
one can recover the original~$H$-structure on~$M$ up to conformal scaling
from the family of submanifolds~$S_y=\lambda\bigl(\rho^{-1}(y)\bigr)$
for~$y\in Y$.

Finally, Merkulov gives representation theoretic criteria on
an irreducibly acting subgroup~$H\subset\GL(\eum)$ with associated 
generalized flag variety~$F\subset\bbP(\eum^*)$ of dimension~$k$ 
and hyperplane bundle~$\cal L$, which guarantee that taking a 
$(2k{+}1)$-dimensional contact manifold~$Y$ and considering the moduli 
space~$M(Y)$ of Legendrian embeddings~$F\subset Y$ that pull back the 
contact bundle~$L$ to be~$\cal L$ yields a smooth moduli space 
endowed with an $H$-structure and a torsion-free connection.

This last step is extremely important, for it provides a way to 
determine which irreducibly acting subgroups~$H\subset\GL(\eum)$
can occur as torsion-free holonomy in terms of representation theory,
specifically, in terms of the vanishing of certain $H$-representations
constructed functorially from~$\eum$.  This provides a 
different approach to solving the torsion-free holonomy problem, one 
that was carried out successfully by a combination of efforts of Chi,
Merkulov, and Schwachh\"ofer.  In particular, this approach led
to the discovery of the remaining groups in Table~IV and, finally,
the proof that Tables I, II, III, and IV exhaust the possibilities
for irreducibly acting torsion-free holonomy.

\vskip 0.5truecm
\noindent
{\bf 3.2.} {\smc Poisson constructions.}  The straightforward
application of exterior differential systems to the holonomy problem
outlined in~\S0.6 does not work for the entries of Table~IV, at 
least at the level described so far.  To see where the problem is,
recall the structure equations derived so far for a torsion-free
connection with holonomy~$H\subset\GL(\eum)$.  They are
$$
\eqalign{
d\omega &= -\theta\w\omega\cr
d\theta &= -{\ts{1\over2}}[\theta,\theta] 
             + {\ts{1\over2}}\,R(\omega\w\omega)\cr
dR &= -\theta.R + R'(\omega)\cr
}
\leqno{(3.1)}
$$
where~$R: B\to K(\euh)$ and ~$R':B\to K^1(\euh)\subset K(\euh)\otimes\eum^*$
are as already defined.

When one considers the first entry of Table~IV, where~$H=\SL(2,\bbR)$
and~$\eum\simeq V_3=\Sym^3(V_1)$, it is not difficult to see that~$K(\euh)
\simeq V_2 = \Sym^2(V_1)$ has dimension~$3$ and that~$K^1(\euh)\simeq V_3$
has dimension~$4$. Its Cartan characters are $(s_1,s_2,s_3,s_4)=(3,1,0,0)$
but, as is easily computed, the 
prolongation~$K^2(\euh)\subset K^1(\euh)\otimes\eum^*$ has 
dimension~$1$, so the tableau is not involutive and Cartan's existence
theorem cannot be applied at this level.  

However, there is a quadratic 
map~$Q:K(\euh)\to K^1(\euh)\otimes\eum^*$ so that the exterior
derivative of the third equation in (3.1) 
is~$\bigl(dR' + \theta.R' - Q(R)(\omega)\bigr)(\omega) = 0$,
implying that there is a function~$R'':B\to K^2(\euh)$
so that the equation
$$
dR' = - \theta.R' + \bigl(R'' + Q(R)\bigr)(\omega)
\leqno{(3.2)}
$$
holds.  Moreover, it is possible to choose the quadratic map~$Q$ in a
unique way so that differentiating this last equation yields $dR''(\omega)=0$.
Now the second prolongation of~$K^1(\euh)$ vanishes, 
so this forces the structure equation
$$
dR'' = 0.
\leqno{(3.3)}
$$
Obviously, differentiating this equation will yield no new information.

At this point, Cartan's general existence theorem for coframings
satisfying prescribed differential identities (a generalization of
Lie's third fundamental theorem) can be applied to the entire
ensemble
$$
\eqalign{
d\omega &= -\theta\w\omega,\cr
d\theta &= -{\ts{1\over2}}[\theta,\theta] 
             + {\ts{1\over2}}\,R(\omega\w\omega),\cr
     dR &= -\theta.R + R'(\omega),\cr
    dR' &= - \theta.R' + \bigl(R'' + Q(R)\bigr)(\omega),\cr
   dR'' &= 0.\cr
}
\leqno{(3.4)}
$$
His theorem implies that for every~$(R_0,R'_0,R''_0)
\in K(\euh){\times}K^1(\euh){\times}K^2(\euh)$, there
is a torsion-free~$H$-structure~$B\subset F(\eum,\eum)$, unique 
up to local diffeomorphism, so that, at a point~$u_0\in B$ one has
$R(u_0)=R_0$, $R'(u_0)=R'_0$, and $R''(u_0)=R''_0$.  Consequently,
the space of diffeomorphism classes of germs of such torsion-free 
$H$-structures is finite dimensional.

What Chi, Merkulov, and Schwachh\"ofer show in~[CMS] is that this exact same 
picture holds for each entry in Table~IV.  Namely, $K(\euh)\simeq\euh$,
$K^1(\euh)\simeq\eum$, ~$K^2(\euh)\simeq\bbF$ ($=\bbR$ or~$\bbC$), 
and there is a quadratic map~$Q:\euh\to\eum\otimes\eum^*$
so that the above analysis of the structure equations goes through
exactly the same as for the original two cases (with the obvious
interpretations of the pairings).

Cartan's theorem then applies and yields not only the existence of 
torsion-free connections with these prescribed holonomies, but that the
space of germs of such $H$-structures, modulo diffeomorphism, is finite 
dimensional.  

In the original two cases that I treated, the finer understanding
of the moduli space of solutions entailed understanding how the
images~$(R,R',R'')(B)\subset\euh{\times}\eum{\times}\bbF$ partition
$\euh{\times}\eum{\times}\bbF$ into subsets.  This analysis would have been 
hopeless were it not for several (at the time) amazing identities that 
I found by brute force calculation.  They even allowed me to prove
existence without using Cartan's existence theorem.  

What is shown in [CMS], however, is that these mysterious identities 
can be explained in terms of a natural Poisson structure on the space
$\euh{\times}\eum{\times}\bbF$  (actually, they regard the last factor
as a parameter and consider Poisson structures on~$\euh{\times}\eum$).
The images~$(R,R',R'')(B)$ turn out to be the
symplectic leaves of this Poisson structure and this point of view
simplifies the reconstruction of the $H$-structure from the
leaf image (though it does not entirely remove some of the global 
difficulties having to do with the symplectic realizations necessary
in their construction).  

\vskip 0.5truecm
\noindent
{\bf 3.3.} {\smc Algebraic classification.}  Once the full list of
the irreducible torsion-free holonomies was known, there arose the 
question of whether this list could be derived through Berger's original 
approach, i.e., representation theory.  Schwachh\"ofer~[Schw1] has shown that
this can indeed be done.  (As is so often the case, knowing the answer 
in advance helps to organize the proof.)  His fully algebraic classification%
\footnote{$^{9}$}{\baselineskip=9pt{\ninerm 
Schwachh\"ofer~[Schw2] informs me that~$\Spin(6,\bbH)\subset\GL(32,\bbR)$
must be added to the exotic symplectic list 
that appears in~[Schw1].  Thus, I have included it in Table IV.}\par}
of the irreducibly acting subgroups~$H\subset\GL(\eum)$ that satisfy
Berger's first criterion and the subset of those that also satisfy Berger's
second criterion still involves quite a bit of case checking, but the
general outline of the argument is clear.

What is particularly intriguing is Ziller's observation (reported in~[Schw1]) 
that this list can be constructed by a simple Ansatz starting from the list 
of Hermitian and quaternionic symmetric spaces.  A direct proof of Ziller's
Ansatz would be highly desirable.

\vskip 0.5truecm
\noindent
{\bf 3.4.} {\smc Two leftover cases.}  As of this writing, each entry
in the four Tables, save two, Entries 3 and 4 in Table~III, 
has been treated in the literature and shown to occur as holonomy, either
by twistor methods or exterior differential systems methods.  Existence
proofs by twistor methods have some difficulty when $\eum$ is a complex 
vector space and the group~$H\subset\GL(\eum)$ is of the 
form~$H=G_\bbC{\cdot}H_s$  where~$H_s$ is the semi-simple part 
and~$G_\bbC\subset\bbC^*$ is a one-dimensional subgroup of~$\bbC^*$, acting
as scalar multiplication on~$\eum$.  The method of exterior differential
systems does not have this problem, but each case does require a separate
treatment.  

In my survey article~[Br3], I left these two entries unsettled in the
case where ~$G_\bbC$ had dimension~$1$ because,
at the time, they did not seem that interesting.  Now that they are the last
unsettled cases, it seems to be a good idea to resolve them, so I will
do that here, though, for lack of space, I will not provide details,
just give the results of the Cartan-K\"ahler analysis.

For the first case, $H = G_\bbC{\cdot}\SL(2,\bbR)\subset\GL(2,\bbC)$, 
one must assume that~$G_\bbC\not\subseteq\bbR^*$, otherwise~$H$ will
not act irreducibly on~$\eum\simeq\bbC^2$.  This leaves a one parameter
family of possibilities~$H_\lambda = \{e^{(i+\lambda)t}\ \vrule\ t\in\bbR\}
\subset\bbC^*$ where~$\lambda$ is any real number.  By conjugation,
one can assume that~$\lambda\ge 0$, so I will do this.  It turns out
that there are two cases:  

If~$\lambda=0$, so that~$H_0=S^1$, it is not difficult to compute that
$K(\euh)\simeq V_4{\oplus}V_2{\oplus}V_0$, while~$K^1(\euh)\simeq 
2V_5{\oplus}2V_3{\oplus}2V_1$ and is involutive, with characters
$(s_1,s_2,s_3,s_4) = (9,9,5,1)$.  Moreover, the torsion is absorbable.
By Cartan's theorem, solutions exist and depend on one function
of four variables.

However, if $\lambda>0$, so that~$H_\lambda\simeq\bbR^+$, one computes
that $K(\euh)\simeq V_4{\oplus}V_2$ while~$K^1(\euh)\simeq2V_5{\oplus}2V_3$ 
and is involutive, with characters
$(s_1,s_2,s_3,s_4) = (8,8,4,0)$.  Moreover, the torsion is absorbable.
By Cartan's theorem, solutions exist and depend on four functions
of three variables.

For the second case, $H = G_\bbC{\cdot}\SU(2)\subset\GL(2,\bbC)$, one
must assume that~$G_\bbC\not\subseteq S^1$, otherwise~$H=\U(2)$ will
preserve a metric on~$\eum\simeq\bbC^2$.  This leaves a one parameter
family of possibilities~$J_\lambda = \{e^{(1+i\lambda)t}\ \vrule\ t\in\bbR\}
\subset\bbC^*$ where~$\lambda$ is any real number.  By conjugation,
one can assume that~$\lambda\ge 0$, though I won't need to do this.  
Here there is only one case:
One computes that $K(\euh)\simeq V_4^\bbR{\oplus}V_2^\bbR\simeq \bbR^8$ 
while~$K^1(\euh)\simeq V_5{\oplus} V_3\simeq \bbC^{10}\simeq\bbR^{20}$. 
The tableau is involutive, with characters
$(s_1,s_2,s_3,s_4) = (8,8,4,0)$.  Moreover, the torsion is absorbable.
By Cartan's theorem, solutions exist and depend on four functions
of three variables.

\vskip 1.5truecm
\noindent
\centerline{\bf REFERENCES}
\vskip 0.5truecm

\item{[Al1]} 
{\smc D. Alekseevskii} - 
{\it Riemannian spaces with unusual holonomy groups}, 
Funct.\ Anal.\ Appl.\ {\bf 2}  (1968), 97--105.

\item{[Al2]} 
{\smc D. Alekseevskii} - 
{\it Classification of quaternionic spaces with a transitive solvable
      group of motions}, 
Math.\ USSR-Izv {\bf 9}  (1975), 297--339.

\item{[AS]} 
{\smc A. Ambrose {\rm and} I. Singer} - 
{\it A theorem on holonomy},
Trans.\ Amer.\ Math.\ Soc.\ {\bf 75} (1953), 428--443.

\item{[BB]} 
{\smc  L. B\'erard-Bergery} - 
{\it On the holonomy of Lorentzian manifolds}, 
Proc.\ Symp.\ Pure Math.\ {\bf 54} (1993), 27--40.

\item{[BI]} 
{\smc  L. B\'erard-Bergery {\rm and} A. Ikemakhen} - 
{\it Sur l'holonomie des vari\'{e}t\'{e}s \break
pseudo-riemanniennes de signature $(n,n)$},
Bull.\ Soc.\ Math.\ France {\bf 125} (1997), 93--114.

\item{[Ber1]}  % MR: 18,149a
{\smc M. Berger},
{\it Sur les groupes d'holonomie des vari\'{e}t\'{e}s \`{a} 
connexion affine et des  vari\'{e}t\'{e}s riemanniennes}, 
Bull.\ Soc.\ Math.\ France {\bf 83} (1955), 279--330.

\item{[Ber2]} 
{\smc M. Berger} -
{\it  Les espaces sym\'etriques noncompacts}, 
Ann.\ Sci.\ \'Ecole Norm.\ Sup.\ {\bf 74} (1957), 85--177.

\item{[Bes]}
{\smc A. Besse} - 
{\sl Einstein Manifolds},
Ergebnisse der Mathematik und ihrer Grenzgebiete, 3. Folge, Band 10,
Springer-Verlag, Berlin, 1987.

\item{[BL]} 
{\smc A. Borel {\rm and} A. Lichnerowicz} - 
{\it Groupes d'holonomie des vari\'et\'es riemann\-iennes}, 
C.R.\ Acad.\ Sci.\ Paris {\bf 234} (1952), 1835--1837.

\item{[BG]}
{\smc R. Brown {\rm and} A. Gray} - 
{\it Riemannian manifolds with holonomy group $\Spin(9)$},
in {\sl Differential Geometry {\rm(}in honor of
Kentaro Yano{\rm)}}, Kinokuniya, Tokyo, 1972,  41--59. 

\item{[Br1]} % MR: 89b:53084
{\smc  R. Bryant} -
{\it Metrics with exceptional holonomy},
Ann.\ Math.\ {\bf 126} (1987), 525--576. 

\item{[Br2]} 
{\smc  R. Bryant} -
{\it Two exotic holonomies in dimension four, path geometries, 
and twistor theory}, Amer.\ Math.\ Soc.\ Proc.\ Symp.\ Pure Math.\  
{\bf 53} (1991), 33--88.

\item{[Br3]} 
{\smc  R. Bryant} -
{\it Classical, exceptional, and exotic holonomies: a status report}, 
in {\sl Actes de la Table Ronde de G\'eom\'etrie Diff\'erentielle en 
l'Honneur de Marcel Berger}, Soc. Math. France, 1996, 93--166.

\item{[BCG]}
{\smc  R. Bryant, {\it et al\/}} -
{\sl Exterior Differential Systems,}
MSRI Series {\bf 18}, Springer-Verlag, 1991.

\item{[BS]} % MR~90i:53055
{\smc R. Bryant {\rm and} S. Salamon} - 
{\it On the construction of some complete metrics with exceptional holonomy}, 
Duke Math.\ J.\ {\bf 58} (1989), 829--850.  

\item{[Ca1]} 
{\smc \'E.\ Cartan} - 
{\sl La G\'eometrie des Espaces de Riemann}, 
M\'emorial des Sciences Mathematiques, Fasc.~IX (1925).  
(Especially, Chapitre~VII, Section~II.)

\item{[Ca2]} 
{\smc \'E.\ Cartan} - 
{\it Les groupes d'holonomie des espaces g\'en\'eralis\'es}, 
Acta.\ Math.\ {\bf 48} (1926), 1--42.

\item{[CMS]} % MR: 98b:53023
{\smc Q.-S. Chi, S. Merkulov, {\rm and} L. Schwachh{\"o}fer} - 
{\it On the existence of infinite series of exotic holonomies.}
Inventiones Math. {\bf 126} (1996), 391--411.

\item{[ChS]} 
{\smc Q.-S. Chi {\rm and} L. Schwachh{\"o}fer} - 
{\it Exotic holonomy on moduli spaces of rational curves}, 
Differential Geom.\ App.\ {\bf 8} (1998), 105--134.

\item{[Gr]} 
{\smc M. Gromov} -
{\it Carnot-Carath\'eodory spaces seen from within},
in {\sl Sub-Riemannian geometry}, Progr. Math.~{\bf 144},
Birkh\"auser, Basel, 1996, 79--323.

\item{[Jo1]} 
{\smc D. Joyce} -
{\it Compact Riemannian $7$-manifolds with 
     holonomy ${\rm G}_2$: {\rm I \& II}}, 
J.\ Differential\ Geom.\ {\bf 43} (1996), 291--328, 329--375.

\item{[Jo2]}   % MR: 97d:53052
{\smc D. Joyce} -
{\it Compact $8$-manifolds with holonomy ${\rm Spin}(7)$}, 
Inventiones Math.~{\bf 123} (1996), 507--552. 

\item{[Jo3]} 
{\smc D. Joyce} -
{\sl Compact manifolds with special holonomy}, 
Oxford University Press, Oxford, 2000 (projected).

\item{[KNa]} 
{\smc S. Kobayashi {\rm and} K. Nagano} -
{\it On filtered Lie algebras and geometric structures, {\rm II}}, 
J.\ Math.\ Mech.\ {\bf 14} (1965), 513--521.

\item{[KNo]} 
{\smc S. Kobayashi {\rm and} K. Nomizu} -
{\sl Foundations of Differential Geometry, Vols.\ $1$ \& $2$}, 
Wiley-Interscience, New York, 1963.

\item{[LM]} 
{\smc H. B. Lawson, Jr. {\rm and} M. L. Michelsohn} -
{\sl Spin geometry}, Princeton Math.\ Series {\bf 38}, 
Princeton Univ.\ Press, Princeton, NJ, 1989.

\item{[Me1]} 
{\smc S. Merkulov} -
{\it  Moduli spaces of compact complex submanifolds 
       of complex fibered manifolds}, 
Math.\ Proc.\ Camb.\ Phil.\ Soc. {\bf 118} (1995), 71--91.

\item{[Me2]} 
{\smc S. Merkulov} -
{\it Geometry of Kodaira moduli spaces}, 
Proc.\ Amer.\ Math.\ Soc.\ {\bf 124} (1996), 1499--1506.

\item{[Me3]} 
{\smc S. Merkulov} -
{\it Existence and geometry of Legendre moduli spaces}, 
Math.\ Z.\ {\bf 226} (1997), 211--265.

\item{[MS1]} 
{\smc S. Merkulov {\rm and} L. Schwachh\"ofer} -
{\it Classification of irreducible holonomies of 
    torsion-free affine connections},
Ann.~Math.~(to appear).

\item{[MS2]} 
{\smc S. Merkulov {\rm and} L. Schwachh\"ofer} -
{\it Twistor solution of the holonomy problem},  
395--402, 
The Geometric Universe, Science, Geometry 
and the work of Roger Penrose, 
S.A. Hugget (ed.), Oxford Univ.~Press, 1998.

\item{[PoS]} 
{\smc Y. S. Poon {\rm and} S. Salamon} -
{\it Quaternionic K\"ahler $8$-manifolds with positive scalar curvature},  
J. Differential Geom.\ {\bf 33} (1991), 363--378.

\item{[Sa]} 
{\smc S. Salamon} -
{\sl Riemannian geometry and holonomy groups}, 
Pitman Research Notes in Math., no. 201, 
Longman Scientific \& Technical, Essex, 1989.

\item{[Scho]} 
{\smc J. Schouten} -
{\sl On the number of degrees of freedom of the geodetically moving
systems}, 
Proc.\ Kon.\ Acad.\ Wet.\ Amsterdam\ {\bf 21} (1918), 607--613.

\item{[Schw1]} 
{\smc L. Schwachh\"ofer} -
{\sl On the classification of holonomy representations}, 
Habilitationsschrift, Universit\"at Leipzig, 1998.

\item{[Schw2]} 
{\smc L. Schwachh\"ofer} -
private communication, 12 October 1999.

\item{[Wi]}
{\smc  B. Wilking} -
{\it On compact Riemannian manifolds with noncompact holonomy groups},
1999 preprint.  \hfil\break 
URL: http://wwwmath.uni-muenster.de/math/inst/sfb/about/publ/wilking.ps

\vskip 1truecm
\hskip 8truecm Robert BRYANT \par
\medskip
\hskip 8truecm Duke University \par
\hskip 8truecm Department of Mathematics\par
\hskip 8truecm P.O. Box 90320\par
\hskip 8truecm DURHAM, NC 27708-0320\par
\hskip 8truecm USA\par
\smallskip
\hskip 8truecm E--mail~: {\tt bryant@math.duke.edu}\par

\bye